\documentclass[12pt]{article}

\usepackage{amsmath, amssymb, latexsym, amscd, amsthm} 

\newtheorem{theorem}{Theorem}[section]
\newtheorem{corollary}[theorem]{Corollary}
\newtheorem{lemma}[theorem]{Lemma}

\newtheorem{statement}[theorem]{Statement}

\theoremstyle{definition}
\newtheorem{definition}[theorem]{Definition}

\def\diam{\operatorname{diam}}
\def\id{\operatorname{id}}
\def\Int{\operatorname{Int}}
\def\Ext{\operatorname{Ext}}

\begin{document}

\title{Cantor sets with 
high-dimensional projections}   
\author{Olga Frolkina\footnote{Supported by Russian Foundation of Basic Research; Grant
No.~19--01--00169.}\\
Chair of General Topology and Geometry\\
Faculty of Mechanics and Mathematics\\
M.V.~Lomonosov Moscow State University\\
Leninskie Gory 1, GSP-1\\
Moscow 119991, Russia\\
E-mail: olga-frolkina@yandex.ru
}

\maketitle

\begin{abstract}
In 1994, J.~Cobb
constructed a tame Cantor set in $\mathbb R^3$
each of whose projections into $2$-planes is
one-dimensional.
We show that 
an Antoine's necklace can serve
as an example of a Cantor set
all of whose projections 
are one-dimensional and connected.
We prove that each Cantor set 
in $\mathbb R^n$, $n\geqslant 3$,
can be moved by a small ambient
isotopy so that the projection of 
the resulting Cantor set 
into each $(n-1)$-plane 
is 
$(n-2)$-dimensional.
We show that if 
$X\subset \mathbb R^n$, $n\geqslant 2$,
is a zero-dimensional compactum 
whose 
projection into some plane $\Pi\subset \mathbb R^n$
with $\dim \Pi \in \{1, 2, n-2, n-1\}$
is zero-dimensional,
then $X$ is tame;
this extends some particular cases of the results of 
D.R.~McMillan, Jr. (1964) and D.G.~Wright, J.J.~Walsh (1982).

We use the technique of defining sequences 
which comes back to 
Louis Antoine.

Keywords:
Euclidean space,
projection,
Cantor set,
embedding,
isotopy,
tame Cantor set,
wild Cantor set,
flat cell,
dimension.

MSC: Primary 54F45; Secondary 57N15, 57M30.
\end{abstract}

\section{Introduction}

Any topological space homeomorphic to 
the standard middle-thirds Cantor set 
$\mathcal C\subset I=[0,1]$ 
is called a Cantor set.
Let $f:\mathcal C\to I$
be a continuous surjection; 
its graph $\Gamma (f) \subset \mathbb R^2$
is a Cantor set whose projection to $Oy$-axis coincides with $I$.
The first description of a 
continuous surjection $f: \mathcal C\to I$
was given in 1884 by G.~Cantor
\cite[p.~255--256]{Cantor};
it maps a point 
$t = \frac{x_1}{3} + \frac{x_2}{3^2} + \ldots 
\in\mathcal C$ 
(where each
$x_i\in \{0;2\}$)
to the point
$f(t)= \frac{x_1}{2}+\frac{x_2}{2^2} + \ldots$ 
of the unit segment.
The remark concerning the projection of $\Gamma (f)$
into the vertical line
can be found in L.~Zoretti's note 
\cite{Zoretti}
who described this as a known fact
(``ce resultat bien connu etant acquis...'').
The aim of Zoretti's note was to
disprove a statement made by F.~Riesz who wrote in \cite[p. 651]{Riesz}
that a projection of a closed
hereditarily disconnected subset of plane
to a straight line is again closed and
hereditarily disconnected
(F.~Riesz attributed this proposition to R.~Baire; this is incorrect since the work of Baire contained a different statement). 
Zoretti also noticed that
there exists
a Cantor set such that its projections to a countable set of
lines contain segments \cite[p.~763]{Zoretti}.
Moreover, he claimed to have described a set all of whose projections 
contain segments \cite[p.~763]{Zoretti} 
(see also \cite{Zoretti-ZB}); unfortunately Zoretti did not give enough details even to prove that his set does not contain arcs.

L.~Antoine
constructed 
a Cantor set in $\mathbb R^2$
whose projections coincide
with those of a regular hexagon
\cite[{\bf 9}, p.272; and fig.2 on p.273]{Antoine-FM}.
Other examples 
of Cantor sets in plane 
all of whose projections are segments can be found in
\cite[p.~124, Example]{Cobb-projections},
\cite[Prop.~1]{Dijkstra-van-Mill}.
Using Antoine's idea, we
show that in $\mathbb R^2$,
each polygon contains a Cantor set
with exactly the same projections (Statement \ref{plane}).

By K.~Borsuk \cite{Borsuk}, for $n\geqslant 2$ there exists a 
Cantor set in $\mathbb R^n$ such that its projection to every hyperplane contains a $(n-1)$-dimensional ball, equivalently, has
dimension $(n-1)$. (For $n=3$,
see an alternative description in
\cite[Thm.~6.2]{Bing57};
another proof of Borsuk's result,
using universal surjectivity property of $\mathcal C$, can be found in \cite[Prop.~3.1]{DGW}.)
As a corollary, 
the projection of the Borsuk set to each
$m$-plane, $m\leqslant n-1$,
has dimension $m$.
We remark that
a set with this property can 
be obtained from any
given Cantor set, using a small
ambient isotopy (Statement \ref{moving-2}).
For $n\geqslant 3$, in $\mathbb R^n$
there is no Cantor set 
all of whose projections to $(n-1)$-planes
are convex bodies
\cite[Thm.~3]{Cobb-projections};
generalizations are obtained in
\cite[Thm.~4.7]{DGW},
\cite[Thm.~1]{BCD}.
In $\mathbb R^3$, we construct a convex body $K$ 
 and a Cantor set $\mathcal A\subset K$
 such that the projections of $K$
 and $\mathcal A$ into each line coincide
(Corollary \ref{coinc}).

J.~Cobb~\cite[Thm.~1]{Cobb-projections} constructed a Cantor set in $\mathbb R^3$ such that its projection to
every $2$-plane is $1$-dimensional.
Cobb asked
\cite[p.~126]{Cobb-projections}:
``Could there be Cantor sets all of whose projections
are connected, or even cells?
...given $n>m>k>0$, does there exist a Cantor set in $\mathbb R^n$ such that
each of its projections
into $m$-planes is exactly $k$-dimensional?''
(Following \cite{BDM}, we call these sets
$(n,m,k)$-sets.)
Examples of
$(n, m, m-1)$- and $(n, n-1, k)$-Cantor sets
are given in \cite[Thm.~1]{Frolkina-proj} and  \cite[Thm.~1]{BDM}, respectively.
The known $(n,m,k)$-Cantor sets
\cite{Borsuk},
\cite[Thm.~6.2]{Bing57},
\cite[Thm.~1]{Cobb-projections},
\cite[Thm.~1]{Frolkina-proj},
\cite[Thm.~1]{BDM}
are tame by constructions and by
\cite[Thm.~I.4.2]{Keldysh} (see Statement~\ref{cell-tame};
for the notions of tame and wild, see
Definition \ref{tame-C}).

Consideration of wild Cantor sets 
provides another approach
to Cobb's question. 
We show that there is
an easily described (wild) Cantor set in $\mathbb R^3$
--- a well-known Antoine's Necklace ---
all of whose plane projections are connected, one-dimensional
(Theorem \ref{basic-example}); and
no of its projections can be ho\-meo\-morphic
to a graph 
(Corollary \ref{graph}).
This answers 
the first part of Cobb's question
for the case of
$3$-dimensional space.

Generalizing Antoine's construction,
W.A.~Blankinship 
\cite{Blankinship}
and
A.A.~Iva\-nov
\cite{Ivanov-diss}
described,
 for $n\geqslant 3$,
(wild) Cantor sets in
$ \mathbb R^n$ with 
non-simply connected complement.
In $\mathbb R^n$, $n\geqslant 3$,
there exist 
wild Cantor sets with simply connected complements
\cite{Kirkor}, 
\cite{deG-O}, \cite{Skora}.
We show that 
each Cantor set in $\mathbb R^n$ can be
slightly moved 
so that the resulting set is an
$(n,n-1,n-2)$-set
(Theorem \ref{moving-1}).
Projections
of a wild Cantor set have rather 
complicated structure
(Statement \ref{projections1}, Corollary \ref{graph}).

In terms of projections, several tameness conditions are known.
W.A.~Blankinship proved that 
if a compactum $K\subset\mathbb R^n$, $n\geqslant 3$,
has a zero-dimensional projection into at least one
hyperplane, then 
$\pi _1 (\mathbb R^n - K) =0$
\cite[Thm.~3E]{Blankinship}.
D.G.~Wright and J.J.~Walsh proved the following:
let $p:\mathbb R^n\to \mathbb R^{n-1}$
be the projection onto a fixed hyperplane;
if $K\subset \mathbb R^n$, $n\geqslant 4$,
is a compactum
with
$\dim K\leqslant n-3$,
$\dim p(K)\leqslant n-3$,
and
$\dim (K\cap p^{-1}(z)) \leqslant 0$
for each $z\in \mathbb R^{n-1}$,
then $K$ is 1-LCC in $\mathbb R^n$;
see \cite[Thm. 5.3, 5.4]{Wright} for $n\geqslant 5$ and \cite{Walsh-Wright} 
for $n=4$.
(By definition, $K\subset \mathbb R^n$ is 1-LCC embedded 
iff for each
$x \in K$ and each neighborhood $U$ of $x$ in $\mathbb R^n$, 
there exists a smaller neighborhood $V$ of $x$ in $\mathbb R^n$
such that any map $\gamma : S^1 \to V - X$ is null-homotopic in $U-X$).
For a zero-dimensional compactum $K\subset \mathbb R^n$ the 1-LCC property is equivalent 
to its tameness in the sense of Definition \ref{tame-C};
see
\cite{Homma} or \cite{Bing57} for $n=3$, 
\cite{McMillan-Taming} for $n\geqslant 5$, 
and
\cite{Shtanko4}, \cite{BDVW} for $n=4$;
for further details,
refer 
to \cite{Edwards} or \cite[3.4]{DV}.
D.R.~McMillan, Jr.
showed that
if a zero-dimensional compactum
$X\subset \mathbb R^{n+m}=\mathbb R^n\oplus \mathbb R^m$, $n,m\geqslant 1$,
satisfies
$\dim p_1 (X) = 0 = \dim p_2(X)$,
where 
$p_1:\mathbb R^n\oplus \mathbb R^m
\to\mathbb R^n$
and 
$p_2:\mathbb R^n\oplus \mathbb R^m
\to\mathbb R^m$
are the projections,
then $X$ is tame \cite[Cor.~2]{McMillan-Taming}.
We show that if $X\subset \mathbb R^n$, $n\geqslant 2$, is
a zero-dimensional compactum,
$\Pi\subset \mathbb R^n$
is a subspace
such that
$\dim p_{\Pi } (X) = 0$ 
and $\dim \Pi \in \{1, 2, n-2, n-1\}$,
then $X$ is tame (Theorem \ref{addition}).

J.~Cobb remarked
that
Cantor sets that raise dimension under all projections
and those that do not are both dense in the Cantor sets in $\mathbb R^n$ \cite[p.~128]{Cobb-projections}.
In a forthcoming paper, we show that
all projections of a typical 
(in the sense of Baire category)
Cantor set are Cantor sets,
partially
answering another question of Cobb.

\subsection*{Notation and Conventions}

An \emph{$m$-plane} in $\mathbb R^n$ is any
$m$-dimensional affine subspace.
An $(n-1)$-plane is also called a \emph{hyperplane}.
We consider only orthogonal projections.
For a plane $\Pi \subset \mathbb R^n$,
the projection map
$\mathbb R^n\to \Pi $
is denoted by $p_{\Pi }$.

The usual Euclidean distance
between points $x$ and $y$
of $\mathbb R^n$ is denoted
by $d(x,y)$.
A \emph{ball} in $\mathbb R^n$
(sometimes we refer to it as an $n$-ball)
is a set of the form
$\{ X \in\mathbb R^n \mid d(X,A) \leqslant r \}$,
 where $A\in\mathbb R^n$ and $r>0$.
An \emph{$n$-cell} is any set homeomorphic to 
a ball in $\mathbb R^n$.

A \emph{circle} in $\mathbb R^3$
is a set obtained by rotating a point
around a straight line.

For a topological manifold-with-boundary
$M$,
denote by 
$\Int M$
and $\partial M$
the interior and the boundary of $M$, correspondingly.

Any topological space
homeomorphic to the usual middle-thirds Cantor set 
is called a Cantor set.
(These are exactly non-empty metric zero-dimensional perfect compacta
\cite[Thm.~12.8]{Moise}.)

All maps are continuous 
unless otherwise specified.

Let $I=[0, 1]$.

An \emph{isotopy} of $\mathbb R^n$
(also called ambient isotopy)
is a level preserving ho\-meo\-morphism
$H:\mathbb R^n \times I\cong \mathbb R^n\times I$
such that $h_0 = \id $,
where $h_t : \mathbb R^n\cong \mathbb R^n$
is defined by $H(x,t) = (h_t (x), t)$.
If, moreover, $H$ is a PL homeomorphism,
then we call it a PL isotopy.
(For PL topology, see e.g. \cite{RS};
an overview of PL notions and results can be found in \cite[1.6]{Rushing}).
We will write $\{ h_t , t\in I \}$ 
or briefly $\{ h_t \}$
instead of $H$.
An \emph{$\varepsilon $-isotopy}
of $\mathbb R^n$
is an isotopy
$\{h_t \}: \mathbb R^n\to\mathbb R^n$
such that
$d(h_t (x), x) \leqslant \varepsilon $
for each $x\in\mathbb R^n$ and each $t\in I$.

A subset $P\subset \mathbb R^n$
is called a \emph{polyhedron}
(or \emph{polyhedral})
if it is the union of a finite collection of simplices.

\section{Basic example: Antoine's Necklace}

The main result of this section is 
Theorem \ref{basic-example}.
I decided to discuss the basic example in
a separate section;
the description and proof are rather
elementary, easier than
in the general case 
(Section \ref{moving-Cantor})
where results from
general theory of 
zero-dimensional compacta are needed.

\begin{definition}\label{tame-C}
A zero-dimensional compact set $K\subset \mathbb R^n$ is called \emph{tame}
if there exists
a homeomorphism $h$ of $\mathbb R^n$ 
onto itself such that 
$h(K)$ is a subset
of a straight line in $\mathbb R^n$;
otherwise, $K$ is called \emph{wild}.
\end{definition}

By Statement \ref{aux2},
we may 
replace
``a homeomorphism $h$ of $\mathbb R^n$'' by ``an isotopy
$\{ h_t\}$ of $\mathbb R^n$'' in Definition \ref{tame-C}; the new definition is
equivalent to the original one.

In $\mathbb R^2$ each zero-dimensional compactum
is tame 
\cite[{\bf 75}, p.~87--89]{Antoine-diss}
(one may also refer to
\cite[Cor.~II.3.2, Cor.~II.3.3]{Keldysh} or
\cite[Chap.~13]{Moise}).

L.~Antoine in \cite{Antoine}
sketched and in
\cite[{\bf 78}, p.~91--92]{Antoine-diss}
explicitly constructed
a Cantor set 
in $\mathbb R^3$ which is now widely known as
an Antoine's necklace;
Antoine proved its wildness
in
\cite{Antoine-diss} (see Statement
\ref{Antoine-separ}).
Let us describe Antoine's construction.

\begin{definition}\label{torus-r-R}
Let $\Pi $ be a $2$-plane in $\mathbb R^3$.
Let $D\subset \Pi $ be a disk of radius $r>0$
with center $Q$, and $\ell \subset \Pi $ a straight line such that $d(Q,\ell ) = R > r$.
A \emph{standard solid torus} $T$
is the solid torus of revolution 
 generated by
revolving $D$ in $\mathbb R^3$ about $\ell $.
The \emph{central circle} of $T$
is the circle generated by rotating the point $Q$. The \emph{center} of $T$ is the center
of its central circle.
\end{definition}

\begin{definition}\label{A-def}
\emph{A simple chain}
in a standard solid torus $T\subset \mathbb R^3$
is a finite family 
$T_1,\ldots , T_q$, $q\geqslant 3$, 
of pairwise disjoint 
congruent
standard solid tori such that

1) $T_1\cup\ldots\cup T_q\subset \Int T$;

2) centers of $T_1,\ldots, T_q$
are subsequent vertices of a regular convex $q$-gon inscribed in the
central circle of $T$;

3) $T_i$ and $T_j$ are linked for
$|i - j |\equiv 1\mod q$, and are not linked
otherwise;

4) for each $i$, the central circle of $T_i$
is zero-homotopic in $T$.
\end{definition}

The chain $T_1,\ldots , T_q$ 
looks like a usual ``necklace'' which
winds once around the hole of $T$;
no one of the tori $T_i$ embraces the hole of $T$.

Antoine takes a standard
solid torus $T$ in $\mathbb R^3$;
he assumes that the ratio $\frac{r}{R}$
is small enough
(where $r$ and $R$ are the numbers described
in Definition \ref{torus-r-R}).
There exists an integer $k$ (sufficiently large 
in comparison with $\frac{r}{R}$) such that
a simple chain of cyclically linked $k$ tori, each 
\emph{geometrically similar} to $T$, can be placed inside $T$,
so that 
their centers lie on the central circle of $T$
and form a regular convex $k$-gon.
Then, he applies a similarity transformation
to place a chain of $k$ tori in the interior of each torus of the previous level, and so on.
Since diameters of the tori tend 
to zero, the intersection 
of sets constructed on all levels
is a Cantor set.

To generalize
this construction, one
allows the tori to be non-standard 
solid tori,
that is, not necessarily tori of revolution;
also, one varies the number of
tori on each stage.
The positions of the tori may also change. 
To be precise, let us give a definition.

\begin{definition}
\emph{An Antoine's Necklace} is a Cantor set in $\mathbb R^3$
which can be obtained as an
intersection $\mathcal A=\bigcap\limits_{i=1}^\infty   M_i$,
where
each $M_i$ 
is the union of a finite number
of pairwise disjoint standard solid tori 
such that:

1) $M_1\subset\mathbb  R^3$ 
is a standard solid torus;

2) $M_{i+1} \subset \Int M_i$ for each $i\geqslant 1$;

3) for each $i\geqslant 1$ and each component $T$ of $M_i$,
the intersection
$M_{i+1} \cap T$ 
is the union of
solid tori,
which form a simple chain in $T$.
\end{definition}

One thus obtains uncountably many
inequivalently embedded 
Antoine's necklaces \cite[Thm.~2]{Sher68}.
Each of them is wild.
L.~Antoine derived the wildness property
for his particular
family of necklaces
from the next statement
proved in
\cite[{\bf 80-82}, p. 93-96]{Antoine-diss};
with minor changes, 
Antoine's reasoning holds for 
each set that satisfies 
Definition \ref{A-def}:

\begin{statement}\label{Antoine-separ}
Let $\mathcal A\subset \mathbb R^3$ be an Antoine's necklace and
$\Sigma\subset \mathbb R^3$
a set homeomorphic to a $2$-sphere.
Suppose that
$\mathcal A$ intersects each
of the (two) connected components
of $\mathbb R^3 - \Sigma $.
Then 
$\Sigma\cap\mathcal A \neq\emptyset $.
\end{statement}

This implies

\begin{statement}\label{proj-connected}
Each orthogonal projection of each Antoine's
necklace is connected.
\end{statement}

Proof.
Let $\mathcal A\subset \mathbb R^3$
be an Antoine's necklace and $\Pi \subset \mathbb R^3$ a $2$-plane.
Recall that $\mathcal A=\bigcap\limits_{i=1}^\infty M_i$,
where each $M_i$
is the union of standard solid tori 
(Definition~\ref{A-def}).

Suppose 
that for some integer
$N$, 
the projection $p_{\Pi }(M_N)$
is not connected.
One may find a simple closed 
curve 
$L\subset \Pi $ 
such that
$L\cap p_{\Pi }(M_N) = \emptyset $,
while 
$p_{\Pi }(M_N)$ intersects
both connected components
$L_{\Int}$ and $L_{\Ext}$
of $\mathbb R^2 - L$
(recall that each component of $M_N$ is a solid torus of revolution).
Take a large number $h\in\mathbb R$ 
such that
$\mathcal A\cap \{ x_3 = \pm h \} = \emptyset $.
The union
$
L \times [-h, h] \cup
L_{\Int } \times \{ \pm h\}
\subset \mathbb R^3
$ is homeomorphic to 
a $2$-sphere; this contradicts
Statement~\ref{Antoine-separ}.

Consequently, each $p_{\Pi }(M_i)$
is connected.
The intersection of a nested sequence of connected, compact sets is connected. Hence
$
p_{\Pi } (A) = 
p_{\Pi } (\cap _{i=1}^\infty M_i )
= 
\cap _{i=1}^\infty p_{\Pi } (M_i )
$
is also connected.

We have proved that the projection into each $2$-plane
is connected;
this easily implies that 
the projection into any straight line
is connected.

We get an immediate corollary
(compare \cite[Thm.~4.7]{DGW}):

\begin{corollary}\label{coinc}
Let $\mathcal A\subset \mathbb R^3$
be an Antoine's necklace.
There exists a convex body $K_{\mathcal A}\subset \mathbb R^3$
such that
$p_{\ell } (\mathcal A) = 
p_{\ell } (K_{\mathcal A}) $
for each line $\ell \subset \mathbb R^3$.
\end{corollary}

Proof.
Let $K_{\mathcal A}$ be the convex hull of $\mathcal A$.
We only need to show that
$K_{\mathcal A}$ is a body
(that is, has a non-empty interior as a subset of 
$\mathbb R^3$).
Suppose not;
then $K_{\mathcal A}$ lies in a $2$-plane.
Recall that each zero-dimensional
compactum in $\mathbb R^2$
is tame \cite[{\bf 75}, p.~87--89]{Antoine-diss};
thus
$\mathcal A$
is tame, a contradiction.

\begin{theorem}\label{basic-example}
There exists an Antoine's Necklace
$\mathcal A$ in $\mathbb R^3$ 
such that 
for each $2$-plane 
$\Pi \subset \mathbb R^3$,
the set 
$p_{\Pi } (\mathcal A) $ is connected and one-dimensional.
\end{theorem}

This fact, once stated, is easily believed.
But it is not as evident as it may seem.
The number of tori increases quickly;
it may happen that they are not sufficiently thin, and the shadows
of the whole family are
``fat'' (compare Statement \ref{moving-2}).
Therefore we need to take care of the tori 
widths at each stage.

\begin{definition}\label{what-is-a-tube}
Let $\ell \subset \mathbb R^3$ be a straight line and $r>0$. 
A \emph{tube}
of radius $r$
with axis $\ell $ is the set
$t(\ell , r) = \{ P\in\mathbb R^3 \mid d(P,\ell ) \leqslant r \}$.
\end{definition}

\begin{definition}\label{what-is-a-strip}
A \emph{strip}
in $\mathbb R^3$ is
a (closed) set between two parallel $2$-planes.
\end{definition}

\begin{definition}
The \emph{width} $w(X)$
of a set $X\subset \mathbb R^n$
is the infimum of the
distances between two parallel hyperplanes
such that 
the part of space between them contains $X$.
\end{definition}

For example, 
$w(t(\ell , r)) = 2r $.

\begin{statement}\label{tubes}
Let $t_1,\ldots , t_s$ 
be a collection of tubes
in $\mathbb R^3$ such that
$w(t_1) + \ldots + w(t_s) < 2\varepsilon $.
Then, for each $2$-plane $\Pi \subset \mathbb R^3$,
the projection
$p_{\Pi } (t_1\cup \ldots \cup t_s)$
contains no $2$-ball of radius $\varepsilon $.
\end{statement}

Proof.
Suppose that  there exists
a $2$-plane $\Pi $ and 
a $2$-ball $B\subset \Pi $ 
of radius $\varepsilon $
such that
$$B\subset p_{\Pi } (t_1\cup \ldots \cup t_s)
= p_{\Pi }(t_1)\cup \ldots \cup p_{\Pi } (t_s).$$
For any tube $t(\ell , r)$,
the projection
$p_{\Pi } (t(\ell , r))$
is either a $2$-ball of radius $r$,
or a strip of width $2r$.
In either case, 
$p_{\Pi } (t(\ell , r))$
can be covered by a strip of width $2r$.
Hence the ball
$B$ of radius $\varepsilon $
can be covered by strips
of total width less than $2\varepsilon $;
this contradicts the
Plank Theorem 
(for the planar case, see
English translation of
A.~Tarski's and H.~Moese's papers (1931--32)
in \cite[Chapter 7]{MMS} where the history
of the problem is also discussed;
the high-dimensional case is treated in \cite{Bang}).

\begin{statement}\label{circle}
Let $L\subset \mathbb R^3$
be a circle.
For each $\varepsilon >0$,
there exists
a finite family of tubes
$t_1,\ldots , t_s$
such that
$L\subset \Int t_1\cup \ldots \cup \Int t_s$
and 
$w(t_1) + \ldots + w(t_s) < \varepsilon $.
\end{statement}

Proof.
Denote the radius of $L$ by $r$.
Let $L_N$ be a regular convex $N$-gon inscribed in $L$. 
Around each straight line containing a side of $L_N$, construct a tube of width 
$4\left( r- r\cos\frac{\pi }{N} \right)$.
It can be easily seen that
the interiors of these tubes cover $L$.
Since the total width
$
N \times 4\left( r- r\cos\frac{\pi }{N} \right)
$
tends to zero as $N\to \infty $,
the result follows.

Proof of Theorem \ref{basic-example}.
Recall that a subset of $\mathbb R^n$
is $n$-dimensional iff
it contains an $n$-ball 
\cite[Theorems 1.8.10, 4.1.5]{Engelking}.
Equivalently,
a subset of $\mathbb R^n$
is $n$-dimensional iff
for some integer $i$ it contains an $n$-ball of radius 
$\frac{1}{i}$.
The desired
Antoine's necklace is constructed in countable number of steps.
On each step $i\in \mathbb N$,
we first
construct 
arbitrary ``preliminary''
simple chains $M_i$ as described
in Definition \ref{A-def};
and then,
we replace them 
(without changing their central circles)
by 
thinner ones so that
no projection of this new $M_i$
contains a $2$-ball of radius $\frac{1}{i}$.
(This is achieved with the help of
Statements \ref{tubes}  and  \ref{circle}.)
Only after that, we pass on to the step $i+1$.
We also require that the diameters 
of the tori which constitute $M_i$
tend to zero as $i\to\infty $.
Then $\mathcal A := \bigcap\limits _{i=1}^\infty M_i$
is an Antoine's necklace in the sense of 
Definition
\ref{A-def}.
The projection
$p_{\Pi}(\mathcal A) = p_{\Pi } \left(\cap _{i=1}^\infty M_i \right)
=\cap _{i=1}^\infty p_{\Pi }( M_i )$
contains no $2$-ball, 
thus $\dim p_{\Pi }(\mathcal A) 
\leqslant 1$.
To finish the proof, apply Statement \ref {proj-connected}.

\section{Additional known facts about Cantor sets in $\mathbb R^n$}

Before we pass on to the higher-dimensional case,
let us discuss some results
which will be of use below.

\begin{definition}
\cite[Def.~I.3.2]{Keldysh}
A set $X$ in an
$n$-dimensional topological manifold $M$
is called
\emph{cellularly separated} (in $M$) if
for each open neighborhood $V$  of $X$
there exists 
a family $\{ u_\alpha \}$ of open subsets of $M$
such that
$X\subset 
\bigcup\limits_{\alpha }  u_\alpha 
\subset 
\bigcup\limits_{\alpha }  \overline u_\alpha  
\subset V$;
each $u_\alpha $ is homeomorphic to $\mathbb R^n$;
each
$\overline u_\alpha $
is an $n$-cell;
and 
$\overline u_\alpha  \cap \overline u_\beta = \emptyset $
for $\alpha \neq \beta $.
\end{definition}

The next proposition was
proved in 
\cite[Thm.~I.4.2]{Keldysh}
(for $n=3$, see 
\cite[Thm.~1.1]{Bing57};
for arbitrary $n$, one may also refer to
\cite[Thm.~1]{Osborne1}).

\begin{statement}\label{cell-tame}
For each $n\geqslant 1$ 
and each zero-dimensional
compact set
$K\subset \mathbb R^n$ the
following conditions are equivalent:

(a) $K$ is cellularly separated in $\mathbb R^n$;

(b) $K$ is tame in $\mathbb R^n$.
\end{statement}

The proof of  \cite[Thm.~I.4.2]{Keldysh} 
implies:

\begin{statement}\label{aux2}
Let $U\subset \mathbb R^n$, $n\geqslant 1$,
be an open connected set,
and $K_1, K_2 \subset U$
be two Cantor sets tame in $\mathbb R^n$.
There exists an isotopy
$\{ h_t \}$ of $\mathbb R^n$ such that
$h_t =\id $ on $\mathbb R^n - U$ for each $t\in I$, and 
$h_1(K_1)=K_2$.
\end{statement}

The following definition is essentially
the one given in 
\cite[p. 662]{Antoine}
(where the term ``vari\'et\'es de d\'efinition de l'ensemble'' is used) and
\cite[pp. 79, 82]{Antoine-diss}
(``surfaces de d\'efinition'').

\begin{definition}
Let $K$ be a zero-dimensional compact 
subset of $\mathbb R^n$. 
A sequence $\{ M_i, i\in\mathbb N\}$ of subsets of $\mathbb R^n$ is called
\emph{a defining sequence for $K$}
if
all
$M_i$'s
are compact 
polyhedral 
$n$-manifolds-with-boundary, 
$M_{i+1} \subset \Int M_i$,
and
$K=\cap_{i=1}^\infty  M_i$.
(The equality $\dim K=0$
implies that
the maximal diameter of the components of $M_i$
necessarily
tends to zero as $i\to\infty $.)
\end{definition}

L.~Antoine
showed that 
each zero-dimensional compactum 
in $\mathbb R^n$, $n\geqslant 1$,
has a defining sequence
\cite[p. 662]{Antoine}, \cite[{\bf 69}, p.78--80]{Antoine-diss}.
Using this result, he proved that
each zero-dimensional compactum 
in $\mathbb R^n$, $n\geqslant 2$,
can be extended to a simple arc
\cite[{\bf 72}, p. 82--84]{Antoine-diss}
(this was first stated in \cite{Riesz}).
In particular, there exists 
a Jordan arc in $\mathbb R^3$
which 
contains an Antoine's necklace; 
this is a first example of a wild 
arc \cite[{\bf 83}, p. 97]{Antoine-diss}.
In \cite{Antoine-173},
Antoine announced
and in \cite{Antoine-FM}
described in detail 
an embedded $2$-sphere in $\mathbb R^3$
which contains 
an Antoine's necklace,
the description can be found also in
\cite[Thm.~18.7]{Moise};
this is the first example of a wild surface.

We will need a stronger property
of defining sequences.
By \cite[Lemma 4]{Armentrout} 
each zero-dimensional compact set $K\subset \mathbb R^3$
has a defining sequence 
$\{ M_i\}$
such that each 
$M_i$
is a disjoint union 
of
polyhedral cells-with-handles.
In general, 
Shtan'ko--Bryant demension theory 
\cite[Prop. 1.2, Thm. 1.4]{Edwards},
\cite[Thm. 3.4.11, 3.4.12]{DV}
implies 

\begin{statement}\label{dem}
Each zero-dimensional compact set $K\subset \mathbb R^n$, $n\geqslant 2$,
has a defining sequence 
$\{ M_i\}$
such that each
$M_i$
is a 
regular neighborhood 
of an $(n-2)$-dimensional
polyhedron;
in particular, 
$M_i$ is a PL manifold-with-boundary.
\end{statement}

Finally, 
let us state one more fact
which we will make use of
(see 
\cite[Cor.~1]{McMillan-Taming}, \cite[Theorem~3]{Osborne1}; 
a very short argument
which
deduces this from the Klee flattening theorem
\cite[Thm.~2.5.1]{Rushing}, \cite[Cor. 2.5.3]{DV}
can be 
found in \cite[Thm.~2]{ZS}).
This proposition is covered by
\cite[Thm. 5.3, 5.4]{Wright}, \cite{Walsh-Wright}.
We obtain another generalization of Statement
\ref{hyperplane} below (see Corollary \ref{Klee-gener}).

\begin{statement}\label{hyperplane}
If a zero-dimensional compactum
$K\subset\mathbb R^n$
lies in a hyperplane,
then $K$ is tame.
\end{statement}

\section{Getting rid of $(n-1)$-dimensional 
projections}\label{moving-Cantor}

\begin{theorem}\label{moving-1}
Let $K\subset \mathbb R^n$, $n\geqslant 2$,
be any Cantor set.
For each $\varepsilon >0$
there exists an $\varepsilon $-isotopy
$\{ h_t \} :\mathbb R^n\cong \mathbb R^n$
such that
$\dim p_{\Pi }(h_1(K)) = n-2$
for each $(n-1)$-plane $\Pi \subset \mathbb R^n$.
\end{theorem}

If $K$ is tame, the result follows
from known constructions.
To prove Theorem \ref{moving-1} for the case of a wild set $K$,
we proceed essentially as we did above,
with Antoine's necklace.
We choose a defining sequence
all of whose elements are regular neighborhoods
of $(n-2)$-dimensional polyhedra;
``compressing''
the elements of the defining sequence
to thin neighborhoods of these subpolyhedra,
we will eliminate $(n-1)$-dimensional projections.
This process 
resembles that in
\cite[Thm. 8]{ALM}.
This gives the inequality 
$\dim p_{\Pi }(h_1(K)) \leqslant n-2$;
the equality will then follow from
\cite{Wright}, \cite{Walsh-Wright}.

Let us extend Definitions \ref{what-is-a-tube}
and \ref{what-is-a-strip}.

\begin{definition}\label{n-tube}
Let $L \subset \mathbb R^n$
be an $(n-2)$-dimensional plane
and $r>0$.
A \emph{tube} of radius $r$ with base $L$
is the set
$t(L,r) = \{ P \in\mathbb R^n \mid
d(P,L) \leqslant r \} $.
\end{definition}

\begin{definition}
An \emph{$n$-strip} 
is a (closed) set between two parallel hyperplanes
in $\mathbb R^n$.
\end{definition}

\begin{statement}\label{tube-projection}
Let $t(L,r)$ be a tube 
of radius $r$ with base $L$ in $\mathbb R^n$,
and let $\Pi \subset \mathbb R^n$
be an $(n-1)$-plane.
The set $p_{\Pi } (t(L,r))$
can be covered by an $(n-1)$-strip in $\Pi $
of width $2r$.
\end{statement}

Proof.
For any point $P\in t(L,r)$ we have
$d(P,L)\leqslant r$.
This implies
${d(p_{\Pi }(P) , p_{\Pi }(L)) \leqslant r}$,
and
$
p_{\Pi }( t(L,r) ) 
\subset \{ Q \in \Pi \mid 
d(Q, p_{\Pi }(L)) \leqslant r \} 
$.
It can be easily seen that the last
set 
is contained in an $(n-1)$-strip in $\Pi $
of width $2r$
(for this, consider two
cases: $\dim p_{\Pi }(L) $ equals $n-2$ or $n-3$).

\begin{statement}\label{shrink}
Let $U\subset \mathbb R^n$, $n\geqslant 2$,
be an open set; let
$M\subset U$ be 
a regular neighborhood
of an $(n-2)$-dimensional polyhedron
$S\subset M$.
Then for each $r >0$
there exists
a PL isotopy
$\{ h_t  \} :\mathbb R^n\cong \mathbb R^n$
such that
$h_t |_{\mathbb R^n - U} = \id $ for each $t\in I$,
$h_1(M)\subset \Int M$,
and $h_1(M)$
can be covered by the interiors of finitely 
many tubes $t(L_1,\delta ),\ldots , t(L_s,\delta )$  with $s\cdot \delta \leqslant r$.
In particular, no projection of
$h_1(M)$ 
into an $(n-1)$-dimensional plane
contains an $(n-1)$-ball of radius $r$.
\end{statement}

Proof.
Cover $S$ by a finite number of
$(n-2)$-planes
$L _1 , \ldots , L _s$.
Take a positive number $\delta < \frac{r}{s}$
such that 
the open $\delta $-neighborhood $O(S,\delta )$
of $S$ 
lies in $\Int M$.
Note that $O(S,\delta )$
is covered by the interiors of the tubes
$t(L_1 , \delta )$,..., $t(L_s,\delta )$
and
$s\cdot \delta < r$.
Let $N \subset O(S,\delta )$ be a regular neighborhood of $S$.
There exists
a PL isotopy $\{ h_t \} : \mathbb R^n\cong \mathbb R^n$
such that 
${h _t}|_ {\mathbb R^n - U } = \id $ for each $t\in I$, and
$h _1 (M) = N$  
\cite[Thm. 1.6.4]{Rushing}.
This is the desired isotopy.
By Statement~\ref{tube-projection}
together with Plank Theorem
\cite{Bang},
no projection of $O(S,\delta )$ 
(hence also of $h_1(M)$)
contains an $(n-1)$-ball of radius $r$.

Now we are ready to prove the main result of this section.

Proof of Theorem \ref{moving-1}.
{\sc Case 1.} $K$ is tame.
This case 
reduces to known results as follows.
There are a finite number
of pairwise disjoint open 
connected sets $U_1,\ldots , U_s \subset U$ such that
$K\subset U_1\cup\ldots\cup U_s$
and $\diam U_j < \varepsilon $ for each $j=1,\ldots , s$.
We may assume that each $K\cap U_j$ is non-empty;
thus $K\cap U_j$ is a tame Cantor set.
For each $j$, take a tame $(n,n-1,n-2)$-Cantor set $K_j\subset U_j$
from 
\cite[Thm.~1]{Frolkina-proj} or \cite[Thm.~1]{BDM}.
By Statement \ref{aux2}
there exists an isotopy
$\{ h _t \} :\mathbb R^n\cong \mathbb R^n$
such that
${h_t} |_{\mathbb R^n - U_1\cup\ldots\cup U_s } = \id $ for each $t\in I$, and 
$h _1 (K\cap U_j ) = K_j $ for each $j=1,\ldots ,s $.
This is the desired map.

{\sc Case 2.} $K$ is wild.
By \cite[Thm.~5.3, 5.4]{Wright} and \cite{Walsh-Wright}
it suffices to obtain inequalities
$\dim p_{\Pi } (h_1(K)) \leqslant n-2$
for each $(n-1)$-plane $\Pi \subset \mathbb R^n$.

For this, we apply Statement \ref{shrink}
infinitely many times.

Let $\{ M_i \}$ be a defining
sequence for $K$ such that 
each $M_i$
is a regular neighborhood of
an $(n-2)$-dimensional polyhedron $K_i$.

Let $ \varepsilon _1 = \varepsilon , \varepsilon _2, \varepsilon _3,\ldots $
be a sequence 
sufficiently fast
decreasing to zero. 
(The exact meaning of this 
will be clarified below.)

There exists an integer $i_1$ such that the diameters of
the components of $M_{i_1}$ do not exceed $\varepsilon _1=\varepsilon $.
Apply Statement \ref{shrink}
taking $U:= \Int M_{i_1}$,
$M:= M_{i_1 + 1}$ and $S:= S_{i_1+1}$;
there exists a PL isotopy 
$\{h^{(1)}_t\}$ of $\mathbb R^n$ 
such that 
${h^{(1)}_t} |_{\mathbb R^n - M_{i_1}} = \id $
for each $t\in I$,
$h_1^{(1)} (M_{i_1+1}) \subset M_{i_1+1}$,
and
no projection of $h^{(1)}_1 (M_{i_1 +1})$
contains a unit $(n-1)$-ball.
Note that $\{h^{(1)}_t \}$ is an $\varepsilon _1$-isotopy.

There exists an integer $i_2 > i_1 $ such that the diameters of
the components of $ h^{(1)} _1( M_{i_2} )$ do not exceed $\varepsilon _2$.
There exists a PL isotopy 
$\{ h^{(2)}_t\}$ of $\mathbb R^n$ 
such that 
${h^{(2)}_t} |_{\mathbb R^n - h^{(1)}_1(M_{i_2})} = \id $
for each $t\in I$,
$h_1^{(2)} ( h_1^{(1)}( M_{i_2+1})) \subset h_1^{(1)} (M_{i_2+1})$,
and
no projection of $h^{(2)} _1 (h^{(1)} _1 (M_{i_2 +1 }))$
contains an $(n-1)$-ball of radius $\frac12$.
Note that
$\{h^{(2)}_t \}$ is an $\varepsilon _2$-isotopy.

Continuing in this way,
for each integer $k$  we find 
an integer $i_k>i_{k-1}$ 
such that the diameters of
the components of $h_1^{(k-1)} \circ h_1^{(k-2)}
\circ\ldots\circ  h^{(1)} _1( M_{i_k} )$ do not exceed $\varepsilon _k$.
There exists 
a PL isotopy 
$\{ h^{(k)}_t \}$ of $\mathbb R^n$ such that 
$$h_t^{(k)} |_{\mathbb R^n - h_1^{(k-1)}\circ\ldots\circ h_1^{(1)} (M_{i_k}) }
=\id \quad\text{ for each }\quad t\in I,$$ 
$$h^{(k)}_1 \circ  \ldots \circ h^{(2)}_1 \circ h^{(1)}_1 (M_{i_k +1} )\subset
h^{(k-1)}_1 \circ  \ldots \circ h^{(2)}_1 \circ h^{(1)}_1 (M_{i_k +1} ),$$
and
no projection of $h^{(k)}_1 \circ  \ldots \circ h^{(2)}_1 \circ h^{(1)}_1 (M_{i_k +1} )$
contains an $(n-1)$-ball of radius $\frac{1}{k}$.
Note that
$\{h^{(k)}_t \}$ is an $\varepsilon _k$-isotopy.

Note that on each step we are free to 
choose $\varepsilon _k$ as small as we wish;
we choose them so that 
the sequence $\{ h^{(k)}_t \circ  \ldots \circ h^{(2)}_t \circ h^{(1)}_t  \}$
converges to an isotopy $\{ h_t\}$
of $\mathbb R^n$, see e.g. 
\cite[Lemma I.4.1]{Keldysh}.

Let us show that $\{ h_t\}$ is the required
isotopy.
By construction,
$h_t |_{\mathbb R^n - M_{i_1}} = \id $
for each $t\in I$; thus $d(h_t(x), x) \leqslant \varepsilon _1 = \varepsilon $ for each $x\in \mathbb R^n$.
No projection of $h_1(K)$
contains an $(n-1)$-ball. In fact,
for each $(n-1)$-plane $\Pi \subset \mathbb R^n$
and each $k\in\mathbb N$
we have by construction
$$
h_1 (K) 
\subset h_1 (M_{i_k +1})
\subset
h^{(k)}_1\circ \ldots \circ h^{(1)}_1 (M_{i_k +1}) ;
$$
hence 
$$
p_{\Pi }(h_1 (K)) \subset 
p_{\Pi } (h^{(k)}_1\circ \ldots \circ h^{(1)}_1 (M_{i_k +1})) ;
$$
the last set does not contain an $(n-1)$-ball
of radius $\frac{1}{k}$.
Therefore $p_{\Pi } (h_1(K))$
does not contain $(n-1)$-balls,
consequently 
$\dim p_{\Pi } (h_1(K)) \leqslant n-2$
by \cite[Theorems 1.8.10, 4.1.5]{Engelking}.

\section{Making all projections have
maximal possible dimension}\label{largest}

In this section, we 
extend results of L.~Antoine and K.~Borsuk.

\begin{statement}\label{plane}
Let $L\subset \mathbb R^2$
be the union of a finite number
of $2$-simplices. There exists
a Cantor set $K\subset \mathbb R^2$
such that 
$p_{\ell }(K)  = p_{\ell }(L)$
for each line $\ell \subset \mathbb R^2$.
\end{statement}

Proof.
The union of a finite number of Cantor sets 
in $\mathbb R^n$
is again a Cantor set 
\cite[Thm. 1.3.1]{Engelking}.
Hence it suffices to consider the case
of a $2$-simplex $L$.
For each $x\in \partial L$ take a hexagon
$H_x$ 
(which is meant to be 
taken together with its interior domain)
affinely equivalent 
to a regular hexagon such that
$x\in H_x\subset L$
and $x$ not a vertex of $H_x$.

The triangle $\partial L$ is compact;
there are finitely many sets
$H_{x_1},\ldots , H_{x_s}$ which cover 
$\partial L$.
The union
$H_{x_1}\cup\ldots \cup H_{x_s}$
form ``an interior collar''
for $\partial L$ in $L$.
For each $i=1,\ldots , s$, 
there exists a Cantor set
$K_i\subset H_{x_i}$
such that 
$p_{\ell } (K_i ) = p_{\ell } ( H_{x_i})$
for each straight line $\ell $
($K_i$ is an affine image of the set 
constructed in \cite[{\bf 9}, p.272; and fig.2 on p.273]{Antoine-FM}).
The union $K=K_1\cup \ldots \cup K_s$
is the desired set.

\begin{statement}\label{moving-2}
Let $K\subset \mathbb R^n$
be a Cantor set, $n\geqslant 2$.
For each 
$\varepsilon >0$
there exists an $\varepsilon $-isotopy 
$\{ h_t \} :\mathbb R^n\cong \mathbb R^n$
such that 
${\dim p_{\Pi }(h_1(K)) = \dim \Pi }$
for each plane $\Pi \subsetneqq \mathbb R^n$.
\end{statement}

Proof.
By 
\cite[Theorems 1.8.10, 4.1.5]{Engelking},
a subset of $\mathbb R^d$
is $d$-di\-men\-sio\-nal iff it contains a $d$-ball;
hence it suffices to get the equality 
$\dim p_{\Pi }(h_1(K)) = n-1$
for each $(n-1)$-plane $\Pi \subset \mathbb R^n$.

Let $\{ M_i\}$ 
be a defining sequence 
for $K$. There exists an integer $N$ such that the diameter of each component of $M_N$
is less than $\varepsilon $.
Enumerate all components 
of $M_N$
by $M_N^{(1)},\ldots , M_N^{(s)}$.
We may assume that for every $i$
the set $K\cap M_N^{(i)}$
is non-empty, hence is
a Cantor set.
Let $K_i \subset K\cap M_N^{(i)}$
be a  Cantor set which
is tame in $\mathbb R^n$ 
(it can be constructed with the help
of Statement~\ref{cell-tame}),
and let $B_i \subset \Int M_N^{(i)}$
be a tame Cantor set 
all of whose projections onto hyperplanes are
$(n-1)$-dimensional \cite{Borsuk}.
By Statement \ref{aux2}
there exists 
an isotopy $\{ h_t \} : \mathbb R^n\cong \mathbb R^n$
such that
$h_t |_{\mathbb R^n - \cup _{i=1}^s M_N^{(i)} } = \id $ for each $t\in I$,
and $h_1(K_i) = B_i$ for each $i=1,\ldots , s$.
This is the desired map.

\section{Tameness of a Cantor set with a zero-di\-men\-sio\-nal projection}

The main result of this section
extends particular cases of
\cite[Cor.~2]{McMillan-Taming}.
On the other side, 
the case $\dim \Pi = n-1$ is covered 
by \cite[Thm. 5.3, 5.4]{Wright} and \cite{Walsh-Wright};
we give a different independent proof. 
We provide unified arguments for three cases $\dim \Pi \in\{ 1,n-2,n-1 \}$;
the case $\dim \Pi = 2$ is easy.

\begin{theorem}\label{addition}
Let $X\subset \mathbb R^n$ be a 
zero-dimensional compact set,
$n\geqslant 2$.
Suppose that 
$\dim p_{\Pi } (X) =0$
for some 
plane $\Pi \subset \mathbb R^n$
whose dimension equals
$1$, $2$, $n-2$ or $n-1$.
Then $X$ is tame.
\end{theorem}

\begin{corollary}\label{Klee-gener}
Let 
$X$
be
a zero-dimensional compactum
such that
$
X\subset 
\mathcal C\times\mathbb R^{n-1}
\subset  
\mathbb R\oplus\mathbb R^{n-1} =
\mathbb R^n$, $n\geqslant 2$.
Then $X$ is tame in $\mathbb R^n$.
\end{corollary}

In our proof of Theorem \ref{addition},
we use
the following lemma inspired by
the Klee trick \cite[Thm. 2.5.1]{Rushing}, \cite[Thm. 2.5.1]{DV}.

\begin{lemma}\label{Lemma}
Let $Q$, $\widehat Q$,
$\widetilde Q$ be $k$-cells in $\mathbb R^k$,
$k\geqslant 1$,
such that
$\widetilde Q \subset \Int \widehat Q 
\subset \widehat Q \subset \Int Q$.
Let $L\subset \mathbb R^{\ell }$, $\ell \geqslant 1$,
be a compact PL
$\ell $-manifold-with-boundary
and $O_L$ its open neighborhood.
Then there exists an isotopy
$\{ F_t \}$  of
$\mathbb R^k \times\mathbb R^\ell $
such that
each $F_t$ is of the form
$
(x,y)\mapsto (f_t (x,y), y)
$ for $x\in \mathbb R^k$, $y\in\mathbb R^\ell $,
$F_t |_{\mathbb R^k \times\mathbb R^\ell
 - Q\times O_L } = \id $ for each $t\in I$,
 and
 $F_1 (\widehat Q \times L)
 \subset \widetilde Q \times L$.
\end{lemma}

Proof of Lemma \ref{Lemma}.
Let $N$ be a regular neighborhood of
$L$ in $\mathbb R^{\ell }$
such that
$N\subset O_L$.
The closure
$\overline{N - L}$
is a PL $\ell $-manifold-with-boundary;
its boundary contains $\partial L$.
By the Collar Neighborhood Theorem,
$\partial L$ has a collar in $\overline{N - L}$.
We may therefore assume that
$\partial L\times [0,1]$ is PL embedded in 
$N$,
and this embedding 
restricted over $\partial L \times \{ 0 \}$
coincides with the restriction of the original inclusion $L\subset O_L$.
Take an isotopy $\{ g_t \}$ of $Q$
such that
$g_0=\id$,
$g_t|_{\partial Q} = \id $ for each $t\in I$,
and
$g_1(\widehat Q)\subset \widetilde Q$.
The desired isotopy $\{ F_t\} $
is defined by the formula
$$
(x,y)\mapsto
\begin{cases}
 \left( g_t(x),y\right) \quad\text{for } (x,y)\in Q\times L; 
\\
\left( g_{t\cdot (1-s)} (x), y \right)
\quad\text{for }
(x,y)\in Q\times \left(\partial L\times \{s\}\right) ;
\\
(x,y)\quad\text{otherwise}.
\end{cases}
$$

Proof of Theorem \ref{addition}.
We may assume that 
$\Pi $ is the coordinate subspace
$ \mathbb R^m \times \{ 0\} ^{n-m} $, $m=\dim \Pi $.
For brevity, we write $p=p_{\Pi }$
and $Y=p(X)$.

The case of a $2$-dimensional plane $\Pi $ 
is easy.
In fact, any zero-dimensional compactum in
$\mathbb R^2$ is tame \cite[{\bf 75}, p.~87--89]{Antoine-diss}, \cite[Cor.~II.3.2, Cor.~II.3.3]{Keldysh}, \cite[Chap.~13]{Moise}. Take a homeomorphism
$h:\mathbb R^2\cong \mathbb R^2$
such that
$h( p(X) ) \subset \{0\}\times \mathbb R^1$.
Define a self-homeomorphism
$H$ of $\mathbb R^n=\mathbb R^2\times\mathbb R^{n-2}$ by
$H= h\times\id_{\mathbb R^{n-2}}$;
we have
$H(X) \subset \{0\}\times\mathbb R^{n-1}$.
By Statement \ref{hyperplane}, $X$ is tame.

For the other cases,
we assume that $X\subset \Int I^n$.
We will show that $X$ satisfies condition
(a) of Statement \ref{cell-tame}.
Fix any $\varepsilon >0$.

\emph{Step 1.}
For each $y\in p(X)$,
the preimage
$X \cap p^{-1} (y) $ is
a zero-dimensional compact subset of 
$\{ y\} \times \Int I^{n-m} $.
There exists
a compact
$(n-m)$-manifold-with-boundary $M_y\subset
\Int I^{n-m}$
such that
$X \cap p^{-1}(y) \subset
\{y\} \times M_y $,
and the diameter of each connected component
of $M_y$
is less than 
$\varepsilon $.
Moreover, we will assume that:
if $m=1$ then each
$M_y$
is a PL manifold-with-boundary;
if $m=n-1$ then each $M_y$ is a closed segment;
and if $m=n-2$ then
each $M_y$ is a flat $2$-cell in $\Int I^2$
(see Definition \ref{flat-cell} and Statement \ref{cell-tame}).

\emph{Step 2.}
We have
$
X \cap 
\left( \{ y\} \times
\left(
I^{n-m} - \Int M_y
\right) \right)
  = \emptyset $
for each $y\in Y$.
Both $X$ and
$
\{ y\} \times
\left(
I^{n-m} - \Int M_y
\right)$
are compact sets;
hence there exists
a PL $m$-cell 
$L_y \subset \Int I^m$
such that 
$y\in \Int L_y$,
$
{X \cap \left( L_y \times
\left(
I^{n-m} - \Int M_y
\right)\right)
 = \emptyset }$,
 and $\diam L_y < \varepsilon $.
We have
$
X\cap p^{-1}(L_y) \subset L_y \times \Int M_y
$.
Note that for each subset
$L\subset L_y$
we have
$
X\cap p^{-1}(L) \subset L \times \Int M_y
$.

\emph{Step 3.}
The collection
$\{ \Int L_y , y\in Y\}$
is a cover of $Y$ by open subsets of $\mathbb R^m$;
take its finite subcover $\{ \Int L_{y_1},\ldots ,
\Int L_{y_q} \}$.
Using the equality
${\dim Y =0}$, we will replace 
$ L_{y_1},\ldots , L_{y_q} $
by smaller subsets
which are pairwise disjoint
(but not necessarily
$m$-cells anymore).
In fact, by compactness of $Y$,
one can find a $\delta >0$
with the property:
for each set $D\subset \mathbb R^m$
with 
$\diam D < \delta $ and
$D\cap Y \neq \emptyset $,
there is an $i\in \{1,\ldots , q\}$ such that
$D\subset \Int L_{y_i}$.
Take 
pairwise disjoint
compact PL 
$m$-manifolds-with-boundary
$N_1,\ldots , N_t\subset \Int I^m$
such that
$Y\subset \Int N_1\cup \ldots \cup \Int N_t$,
$\diam N_j < \delta $ 
and $Y\cap N_j \neq\emptyset $
for each $j\in\{1,\ldots , t\}$ (Statement \ref{dem}).
Each $N_j$ is a subset of some 
$\Int L_{y_i}$ with $i\in\{ 1,\ldots , q\}$.

Replace the family $ \{ L_{y_1} , \ldots , L_{y_q} \}$
by $\{ N_1,\ldots , N_t \}$.
Saving our old notation, we will assume
that $L_{y_1} , \ldots , L_{y_q}$
are
themselves 
pairwise disjoint 
compact PL
$m$-manifolds-with-boundary (not 
ne\-ces\-sa\-ri\-ly
$m$-cells).
Take pairwise disjoint 
open neighborhoods
$O(L_{y_1}), \ldots , O(L_{y_q})$
of $L_{y_1},\ldots , L_{y_q}$
in $\mathbb R^{m}$
with $\diam O(L_{y_i}) < \varepsilon $ 
for each $i$.

\emph{Step 4.}
For each $i\in\{1,\ldots , q\}$,
denote 
the connected components of 
$M_{y_i}$
by $M_{y_i,1}, \ldots , M_{y_i , \alpha (i)}$;
we get
$$
X = \bigcup\limits_{i=1}^q
\left(
X\cap p^{-1} (L_{y_i})
\right)
\subset
\bigcup\limits_{i=1}^q
\left(
\Int L_{y_i}
\times
\Int M_{y_i} 
\right)
=
\bigcup\limits_{i=1}^q
\bigcup\limits_{j=1}^{\alpha (i) }
\left(
\Int L_{y_i}
\times
\Int M_{y_i , j}
\right) .
$$

\emph{Step 5.}
We are now ready to accomplish the proof.

\emph{Case 1: $m=1$.}
We may assume that $L_{y_1}, \ldots , L_{y_q}$
are pairwise disjoint closed segments.
For each $i\in \{1,\ldots , q\}$ take closed segments
$\widehat L_{y_i}$ and 
$\widetilde L_{y_i}$
such that
$\widetilde L_{y_i}
\subset \Int \widehat L_{y_i}
\subset \widehat L_{y_i}
\subset \Int L_{y_i}$
and
$X\subset \bigcup\limits_{i=1}^q
\bigcup\limits_{j=1}^{\alpha (i) }
\left(
\Int \widehat L_{y_i}
\times
\Int M_{y_i , j} \right) $.
For each $j \in\{1,\ldots , \alpha (i) \}$
let $O(M_{y_i ,1})$,..., $O(M_{y_i , \alpha (i)})$ be pairwise disjoint 
open neighborhoods of 
$M_{y_i ,1}$,..., $M_{y_i , \alpha (i)}$
in $\mathbb R^{n-m}$,
each of diameter less than~$\varepsilon $.
For each $M_{y_i,j}$ take an $(n-m)$-ball
$B_{y_i , j} \subset \mathbb R^{n-m}$
of radius $\varepsilon \sqrt {n-m}$
such that
$M_{y_i,j} \subset \Int B_{y_i,j}$
(balls $B_{y_i , j}$ may intersect each other).
Using Lemma~\ref{Lemma},
for each $i\in\{1,\ldots , q\}$
and each $j\in \{1,\ldots , \alpha (i)\}$
construct an
isotopy
$\{ F_{(i,j),t} \}$ of $\mathbb R^n$
identical outside
$L_{y_i} \times O(M_{y_i , j})$
such that
$F_{(i,j), 1} (\widehat L_{y_i } 
\times M_{y_i , j})
\subset \widetilde L_{y_i} \times M_{y_i , j}
\subset \widetilde L_{y_i} \times \Int B_{y_i , j}$.
All isotopies $\{ F_{(i,j), t} \}$
``glued together'' give an isotopy $\{ F_t \}$.
The set $\widetilde L_{y_i} \times B_{y_i , j}$
is an $n$-cell,
and
$\diam (\widetilde L_{y_i} \times B_{y_i , j}) <
\sqrt{\varepsilon ^2 + 4\varepsilon ^2 (n-m)}
= \varepsilon \sqrt{1+4(n-m)}$.
Note that
$\diam (L_{y_i} \times O(M_{y_i , j})) <
\varepsilon \sqrt 2$,
hence
$\{ F_{t} \}$ 
is an $\varepsilon \sqrt 2$-isotopy.
Thus 
$\left\{ \left( F_{t}  \right) ^{-1} \right\}$ 
is also an $\varepsilon \sqrt 2$-isotopy,
and
$ \left( F_{1} \right) ^{-1} 
( \widetilde L_{y_i} \times B_{y_i , j} ) $
is an $n$-cell of diameter
less than
$2\sqrt2 \varepsilon +  \varepsilon \sqrt{1+4(n-m)}$.
The cells 
$\left( F_{1} \right) ^{-1} 
(  \widetilde L_{y_i} \times  B_{y_i , j} )$
are pairwise disjoint,
and 
$X\subset 
\bigcup\limits_{i=1}^q
\bigcup\limits_{j=1}^{\alpha (i) }
 \left( F_{1} \right) ^{-1} 
( \Int \widetilde L_{y_i} \times \Int B_{y_i , j} ) $.
Statement~\ref{cell-tame} now
implies that $X$ is tame.

\emph{Case 2: $m=n-1$.}
For each $i\in\{1,\ldots , q\}$
and each $j\in \{1,\ldots , \alpha (i)\}$
take closed segments
$\widetilde M_{y_i , j}$
and $\widehat M_{y_i,j}$ such that
$\widetilde M_{y_i , j} \subset
\Int \widehat M_{y_i,j}
\subset
\widehat M_{y_i,j}
\subset \Int M_{y_i,j} $,
$X\subset \bigcup\limits_{i=1}^q
\bigcup\limits_{j=1}^{\alpha (i) }
\left(
\Int L_{y_i}
\times
\Int \widehat M_{y_i , j} \right) $,
and all $\widetilde M_{y_i , j}$
are pairwise disjoint.
For each $i\in \{1,\ldots , q\}$
take an open neighborhood
$O(L_{y_i})$ of $L_{y_i}$ in $\mathbb R^m$
such that $\diam O(L_{y_i}) < \varepsilon $
and all $O(L_{y_i})$ are pairwise disjoint;
take an $m$-ball $B_{y_i}$
of radius $\varepsilon \sqrt{m}$ 
such that $L_{y_i} \subset \Int B_{y_i}$
(these balls may intersect each other).
Using Lemma \ref{Lemma},
for each $i\in\{1,\ldots , q\}$
and $j\in \{1,\ldots , \alpha (i)\}$
construct an
isotopy
$\{ F_{(i,j),t} \}$ of $\mathbb R^n$
identical outside
$O(L_{y_i}) \times M_{y_i , j}$
such that
$F_{(i,j), 1} (L_{y_i } 
\times \widehat  M_{y_i , j})
\subset  L_{y_i} \times \widetilde M_{y_i , j}
\subset (\Int B_{y_i} )\times \widetilde M_{y_i , j}$.
As in Case 1,
all isotopies $\{ F_{(i,j), t} \}$ glue
together to an isotopy $\{F_t\}$;
the sets 
$\left( F_{1} \right) ^{-1}
( B_{y_i} \times \widetilde M_{y_i , j} )$
are pairwise disjoint $n$-cells
of diameter $<k\varepsilon $ which cover $X$
($k$ depends on $n$ and $m$, and not on $\varepsilon $),
and Statement~\ref{cell-tame} applies.

\emph{Case 3: $m=n-2$.}
We use the same argument as in Case 2.
The only difference is that
$\widetilde M_{y_i , j}$
and $\widehat M_{y_i,j}$ 
are $2$-cells (since any zero-dimensional compactum in
$\mathbb R^2$ is tame);
we also assume that each $\widetilde M_{y_i ,j}$
is flat in $\mathbb R^2$
(see Definition \ref{flat-cell}).
The details are omitted.

\section{Tameness of a Cantor set with ``simple'' projection}

In this last section, we make an attempt to 
investigate
relationship between wildness of
a Cantor set and ``complexity'' of its projections.

\begin{definition}\label{flat-cell}
A $k$-cell $X\subset \mathbb R^n$ 
is called \emph{flat}
if there exists a homeomorphism $h$
of $\mathbb R^n$ onto itself such that
$h(X)$ is a $k$-simplex.
\end{definition}

For more information on flat and tame embeddings,
in particular 
for examples 
of cells which are not flat,
refer to the books 
\cite{Moise}, 
\cite{Keldysh}, \cite{Rushing}, 
\cite{DV}.

\begin{statement}\label{projections1}
Let $K\subset \mathbb R^n$, $n\geqslant 3$ be a zero-dimensional compact
set.
Suppose that for some $m$-plane $\Pi $, 
where $1\leqslant m\leqslant n$,
there
exists a countable family
of subsets $X_1,X_2,\ldots \subset \Pi $
such that
each $X_i$ is a flat cell in $\Pi $,
$\dim X _i \leqslant m-1$,
and $p_{\Pi }(K)\subset X_1\cup X_2\cup \ldots $.
Then $K$ is tame in $\mathbb R^n$.
\end{statement}

Proof.
Let 
$\Pi ^{\bot }$ be the orthogonal complement
of $\Pi $ in $\mathbb R^n$.
Since $K$ is compact,
there exists an $(n-m)$-ball $B\subset \Pi ^{\bot }$
such that
$K\subset \Pi \times B$.
We have 
$
K\subset 
(X_1\times B) \cup 
(X_2\times B) \cup\ldots
$.
We may assume that 
for each $i$ the set $K\cap (X_i\times B)$
is non-empty, hence it is 
a zero-dimensional compactum.
It is easy to see that each $X_i \times B$
is a flat cell in $\mathbb R^n$,
and $\dim (X_i\times B) \leqslant
n-1$. Therefore
$K_i$ is tame in $\mathbb R^n$
(Statement \ref{hyperplane}).
Now $K$
is tame in $\mathbb R^n$
by \cite[Thm.~8]{Osborne1} (for $n=3$ 
see also \cite[Thm.~6.1]{Bing57}).

\begin{corollary}\label{graph}
Let $K\subset \mathbb R^n$ be 
a 
zero-dimensional compact
set, $n\geqslant 3$.
Suppose that for some $2$-plane $\Pi $,
there exists 
a subset $X\subset \Pi $
homeomorphic to a one-dimensional polyhedron
such that
$p_{\Pi }(K) \subset X$.
Then $K$ is tame in $\mathbb R^n$.
\end{corollary}

Proof.
A set homeomorphic to a 
one-dimensional polyhedron 
is the union of a finite family
of simple arcs.
Each simple arc in plane is flat
by Antoine's theorem
\cite[p.21, {\bf 17}]{Antoine-diss}, \cite[Thm. II.4.3]{Keldysh}, \cite[Thm.~10.8]{Moise}.
Statement \ref{projections1}
now applies.

\end{document}